\documentclass[aos,MSNbibl,nameyear,dvips]{arximspdf}
\usepackage{graphicx}

\doi{10.1214/11-AOS891}
\volume{39}
\issue{3}
\pubyear{2011}
\firstpage{1803}
\lastpage{1826}

\makeatletter

\newtheorem{theorem}{Theorem}
\newtheorem{lemma}{Lemma}
\newtheorem{prop}{Proposition}
\newproclaim{remark}{Remark}
\newtheorem{coro}{Corollary}

\makeatother

\begin{document}
\begin{frontmatter}

\title{Estimation of extreme risk regions under multivariate regular variation}
\runtitle{Extreme risk regions}

\begin{aug}
\author[A]{\fnms{Juan-Juan} \snm{Cai}\ead[label=e1]{j.cai@tilburguniversity.edu}},
\author[A]{\fnms{John H. J.} \snm{Einmahl}\corref{}\ead[label=e2]{j.h.j.einmahl@tilburguniversity.edu}}
and
\author[A]{\fnms{Laurens} \snm{de Haan}\ead[label=e3]{ldehaan@ese.eur.nl}\thanksref{aut1}}
\runauthor{J.-J. Cai, J. H. J. Einmahl and L. de Haan}
\thankstext{aut1}{Supported in part by FCT/PTDC/MAT/112770/2009.}
\affiliation{Tilburg University, Tilburg University and Lisbon University}
\address[A]{Department of Econometrics \& OR and CentER\\
Tilburg University \\
P.O.\ Box 90153\\
5000 LE Tilburg\\
The Netherlands\\
\printead{e1}\\
\hphantom{E-mail: }\printead*{e2}\\
\hphantom{E-mail: }\printead*{e3}} 
\end{aug}

\received{\smonth{7} \syear{2010}}
\revised{\smonth{3} \syear{2011}}

%
\begin{abstract}
When considering $d$
possibly dependent random variables, one is often interested in extreme risk
regions, with very small probability~$p$. We consider
risk regions of the form $\{\mathbf{z} \in\mathbb{R}^d \dvtx
f(\mathbf{z})\leq\beta\}$, where $f$ is the joint density
and $\beta$ a small number. 
Estimation of such an extreme risk region is difficult since it
contains hardly any or no data. Using extreme value theory, we
construct a natural estimator of an extreme risk region and prove a
refined form of consistency, given a random sample of multivariate
regularly varying random vectors. In a detailed simulation and
comparison study, the good performance of the procedure is
demonstrated. We also apply our estimator to financial data.
\end{abstract}

%

\end{frontmatter}

\section{Introduction}\label{sec1}

A two-dimensional normal density or Student $t$-density is constant on
boundaries of certain ellipses. Outside such an ellipse the density is
lower than inside. It is straightforward to find such an outer region
and its contour (line), for a given small probability. We can consider
such contour as a natural multidimensional extension of a
(one-dimensional) quantile. Even for extreme sets, that is,\ very low
density levels, the calculations are straightforward.

In this paper we consider, much more general, \textit{multi}variate
regularly varying distributions [for a review, see \citet{JesMik06}]. We consider the latter distributions,
since we want to explore in particular extreme sets, that is,\ sets far
removed from the origin.
A random vector $\mathbf{X}$ is multivariate regularly varying if
there exist
a constant $\alpha>0$, the index and an arbitrary probability measure
$\Psi$ on $\Theta=\{\mathbf{z}\in\mathbb{R}^d\dvtx \| \mathbf{z}\|
=1\}$, the unit
hypersphere, such that
\begin{equation}\label{mrv}
\lim_{t\rightarrow\infty }\frac{\mathbb{P} ( \|\mathbf{X}\|\geq
tx , \mathbf{X}/\|\mathbf{X}
\|\in A)}{\mathbb{P} (\|\mathbf{X}\| \geq t)}=x^{-\alpha}\Psi(A)
\end{equation}
for every $x>0$ and Borel set $A$ in $\Theta$ with $\Psi(\partial
A)=0$, with $ \|\mathbf{X}\|$ the $L_2$-norm of~$\mathbf{X}$; see
Rva\v{c}eva
(\citeyear{Rva62}).
An equivalent statement is
\begin{equation}\label{rad}
\lim_{t\rightarrow\infty }\frac{\mathbb{P} ( \|\mathbf{X}\|\geq
tx )}{\mathbb
{P} (\|\mathbf{X}\| \geq t)}=x^{-\alpha}\qquad  \mbox{for } x>0,
\end{equation}
and there exists a measure $\nu$ such that
\begin{eqnarray}\label{(1.1)}
\lim_{t\rightarrow\infty }\frac{\mathbb{P} ( \mathbf{X}\in
tB)}{\mathbb{P}
(\|\mathbf{X}\| \geq t)}=\nu(B)<\infty
\end{eqnarray}
for every Borel set $B$ on $\mathbb{R}^d$ that is bounded away from
the origin and satisfies $\nu(\partial B)=0$; here $tB=\{t\mathbf
{z}\dvtx\mathbf{z}\in
B\}$. 
Note that $\nu$ is homogeneous,
that is,\ for all
$a>0$,
\begin{eqnarray}
\nu(aB)=a^{-\alpha }\nu(B).\label{(1.2)}
\end{eqnarray}
Clearly, on $\{\mathbf{z}\in\mathbb{R}^d\dvtx \|\mathbf{z}\|\geq1\}$, $\nu$ is a \textit
{probability} measure.
The limit relation in (\ref{(1.1)}) is a multivariate analogue of the
``peaks-over-threshold'' or ``generalized Pareto limit'' method in
one-dimensional extreme value theory. Particular cases of (\ref{mrv})
are distributions in the sum-domain of attraction of $\alpha$-stable
distributions and heavy tailed elliptical distributions such as
multivariate $t$-distributions [see \citet{Has06}].

We require the convergence in (\ref{rad}) and (\ref{(1.1)}) at the density
level:  
\begin{longlist}[(a)]
\item[(a)] Suppose that the distribution of $\mathbf{X}$ has a continuous
and positive density $f$ and that for some positive function $q$ and
some positive function $V$ regularly varying at infinity with negative
index $-\alpha $, we have
\begin{equation}\label{pw}
\lim_{t\rightarrow\infty }\frac{f(t\mathbf
{z})}{t^{-d}V(t)}=q(\mathbf{z})\qquad  \mbox{for
all } \mathbf{z}\neq\mathbf{0}
\end{equation}
and
\begin{equation}\label{(1.3)}
\lim_{t\rightarrow\infty }\sup_{\mathbf{z}\in\Theta} \biggl|\frac
{f(t\mathbf{z}
)}{t^{-d}V(t)}-q(\mathbf{z}) \biggr|=0.
\end{equation}
%
Then $q$ is continuous on $\mathbb{R}^d\setminus\{\mathbf{0}\}$ and
$q(a\mathbf{z})=a^{-d-\alpha }q(\mathbf{z})$ for all $a>0$ and
$\mathbf{z}\neq\mathbf{0}$.
Throughout, \textit{we can and will take} $V(t)=\mathbb{P}(\|\mathbf{X}
\|>t)$ (see Lemma~\ref{lem1}, Section~\ref{sec5}).
From Lemma~\ref{lem1}, it follows that doing so (\ref{(1.1)}) holds with $
\nu(B)=\int_B q(\mathbf{z})\,d\mathbf{z}
$.
\end{longlist}

The extreme region will be of the form
\[
Q=\{\mathbf{z}\in\mathbb{R}^d\dvtx  f(\mathbf{z})\leq\beta \},
\]
where $f$ is the probability density of the random vector $\mathbf
{X}$; $\beta
$ is determined in such a way that the probability of $Q$ is equal to a
given very small number $p$, like $1/10$,$000$.

It is the purpose of this paper to estimate $Q$ based on $n$ i.i.d.\
copies of $\mathbf{X}$. Note that the shape of $Q$ is
not predetermined, it depends on the density~$f$. For the estimation of
$Q$, we will use an approximation of $f$ based on the density of $\Psi
$. The values of $p$ we consider are typically of order $1/n$. This
means that the number of data points that fall in $Q$ is small and can
even be zero, that is,\ we are extrapolating outside the sample. This
lack of relevant data points makes estimation difficult. The estimation
of $Q$ is a multivariate analogue of the estimation of extreme
quantiles in the univariate setting; see, for example,\ de Haan and
Ferreira (\citeyear{deHFer06}), Chapter~4. The multivariate case
is much more complicated, however, since we have to estimate a whole
set instead of only one value.

Having an estimate of $Q$ can be important in various settings. It can
be used as an alarm system in risk management: if a new observation
falls in the estimated $Q$ it is a signal of extreme risk. See
\citet{EinLiLiu09} for an application to aviation safety along
these lines.
In a financial or insurance setting, points on the boundary of the
estimate of $Q$ can be used for stress testing. The estimate of $Q$ can
also be used to rank extreme observations (see Remark~\ref{rem3},
Section~\ref{sec2}).

For the ``central'' part of the distribution, that is,\ $\beta$ is
fixed (and ``not too small''), nonparametric estimation of density
level sets has been studied in depth in the literature. Two approaches
are used, the plug-in approach using density estimation [see
\citet{BalCueCue01} and \citet{RigVer09}], and the excess
mass approach [see M\"{u}ller and Sawitzki (\citeyear{MulSaw91}),
\citet{Pol95} and \citet{Tsy97}]. Our estimation problem
and (hence) our approach are quite different from these.

This paper is organized as follows. In Section~\ref{sec2}, we derive our
estimator and show a refined form of consistency. A simulation and
comparison study is presented in Section~\ref{sec3} and a financial application
is given in Section~\ref{sec4}. Section~\ref{sec5} contains the proof of the main result.

\section{Main result}\label{sec2}

Consider a random sample $\mathbf{X}_1, \mathbf{X}_2,\ldots, \mathbf
{X}_n $ with
$\mathbf{X}_i\stackrel{d}{=}\mathbf{X}$, for $i=1,\ldots, n$; their common
probability measure on $\mathbb{R}^d$ is denoted with $P$. Write $R_i$
for the
radius $\|\mathbf{X}_i\|$ and $\mathbf{W}_i$ for the direction
$\mathbf{X}_i/\|\mathbf{X}_i\|$ of
$\mathbf{X}_i$. We wish to
estimate an extreme risk region of the form
\[
Q=\{\mathbf{z}\in\mathbb{R}^d\dvtx  f(\mathbf{z})\leq\beta \},
\]
where $\beta $ is such that $PQ=p>0$, where $p=p_n\to0$, as $n \to
\infty$.
This means that $Q$ and $\beta $ depend on $n$, that is, $Q=Q_n$ and
$\beta =\beta _n$. 
We shall
connect $Q_n$ to a fixed set $S$ not depending on $n$, defined by
\[
S=\{\mathbf{z}\dvtx q(\mathbf{z})\leq1\}.
\]
It will
turn out that $Q_n$ can be approximated by a properly inflated version
of $S$. In fact, it follows from (\ref{(1.3)}) that the risk regions
are asymptotically homothetic as a function of $p$, for small values of $p$.
Define $H(s)=1-V(s)=\mathbb{P}(R\leq s)$
and
$U(t)=H^{-1}(1-\frac{1}{t})$. Note that $U$ is regularly varying at
infinity with index $1/\alpha $.

We will approximate $Q_n$ in two steps by a (deterministic) region
$\widetilde{Q}_n$. This approximation satisfies
\begin{equation}\label{appro}
\frac{P(Q_n \triangle\widetilde{Q}_n)}{p}\to0
\end{equation}
($\triangle$ denotes ``symmetric difference'') and is based on the
above limit relations. The region $\widetilde{Q}_n$ can therefore be
estimated using extreme value theory.
The first step is to establish an approximation of $\beta =\beta
(p)$. Let
\begin{longlist}[(b)]
\item[(b)] $k=k_n(<n)$ be a sequence of positive
integers such that $k\rightarrow\infty $ and $k/ n\rightarrow
0$.
\end{longlist}

The region $Q_n$ is approximated by
\[
\bar{Q}_n= \biggl\{\mathbf{z}\dvtx f(\mathbf{z})\leq\biggl(\frac{np}{k\nu(S)}
\biggr)^{1+{d}/{\alpha
}}\frac{1}{({n}/{k})(U({n}/{k}))^d} \biggr\}.
\]
Next, we approximate $\bar{Q}_n$ by a further region $\widetilde
{Q}_n$ defined in
terms of the limit density $q$ rather than $f$:
\begin{equation}\label{til}
\widetilde{Q}_n=U\biggl(\frac{n}{k}\biggr)\biggl(\frac{k\nu(S)}{np}\biggr)^{1/\alpha }S.
\end{equation}
Indeed, $S$ and this approximation of $Q_n$ are homothetic.

Write
\[
B_{r,A}=\{\mathbf{z}\dvtx  \|\mathbf{z}\|\geq r, {\mathbf{z}}/{\|\mathbf
{z}\|} \in A\}
\]
for a Borel set $A$ on $\Theta$. Clearly,
$B_{r,A}=rB_{1,A}$ and hence
$
\nu(B_{r,A})=r^{-\alpha }\nu(B_{1,A}).
$
The relation between the spectral measure $\Psi$ and $\nu$ is [cf.\
(\ref{mrv}) and (\ref{(1.1)})]
\[
\Psi(A)=\nu(B_{1,A})
\]
for
a Borel set $A\subset\Theta$. Recall that the spectral measure is a
probability measure.
The existence of a density $q$ of $\nu$ implies the existence of a density
$\psi$ of $\Psi$, that is,
\[
\Psi(A)=\int_A \psi(\mathbf{w})\,d\lambda(\mathbf{w}),
\]
where $\lambda$ is the Hausdorff measure (surface area) on
$\Theta$ and
\[
q(r\mathbf{w})=\alpha r^{-\alpha -d}\psi(\mathbf{w}).
\]
Next, we write $S$ and $\nu(S)$ in terms of the spectral
density:
\[
S= \bigl\{\mathbf{z}=r\mathbf{w}\dvtx  r\geq(\alpha \psi(\mathbf{w}))^{{1}/{(\alpha
+d)}}, \mathbf{w}\in
\Theta\bigr\}
\]
and hence
\[
\nu(S)=\alpha ^{-{\alpha }/{(\alpha +d)}}\int_{\Theta}
(\psi(\mathbf{w}))^{{d}/{(\alpha +d)}}\,d\lambda(\mathbf{w}).
\]

To estimate $\widetilde Q_n$, we need estimators for $U(n/{k}), \alpha
, S$
and $\nu(S)$. From the above expressions
for $S$ and $\nu(S)$, we see that this means that we have to estimate
$U(n/{k})$,
$\alpha $ and $\psi$. First, we define
\[
\widehat{U}\biggl(\frac{n}{k}\biggr)=R_{n-k: n}
\]
[the $(n-k)$th order statistic of the $R_i$, $i=1, \ldots, n$].
Since the tail of the distribution function of
$R$ is regularly varying with index $-\alpha $, we can use one of
the well-known estimators of the extreme value index $1/\alpha $,
based on
the $R_i, i=1, \ldots,
n$; see, for example, \citet{Hil75}, \citet{Smi87} and
\citet{DekEindeH89}. It remains to estimate $\psi$.
Let $K\dvtx [0,1]\rightarrow[0,1]$ be a continuous and
nonincreasing (kernel) function with $K(0)=1$ and $K(1)=0$.
For
$\mathbf{w}\in\Theta$, define an estimator of $\psi(\mathbf{w})$ by
\[
\widehat{\psi}_n(\mathbf{w})=\frac{c(h, K)}{k}\sum_{i=1}^n K\biggl(\frac
{1-\mathbf{w}^T\mathbf{W}_i}{h}\biggr)1_{[R_i>R_{n-k: n}]}
\]
with $0<h<1$ and
\[
c(h, K)=\biggl(\int_{C_\mathbf{w}(h)}K\biggl(\frac{1-\mathbf{v}^T\mathbf{w}}{h}\biggr)\,d\lambda
(\mathbf{v})\biggr)^{-1},\quad
C_\mathbf{w}(h)=\{\mathbf{v}\in\Theta\dvtx \mathbf{w}^T\mathbf{v}\geq
1-h\};
\]
cf.\
\citet{HalWatCab87}.

For estimating $Q_n$ it suffices to estimate $\widetilde{Q}_n$, see
(\ref
{appro}). Hence, in view of~(\ref{til}), we define
\begin{eqnarray}\label{defqhat}
\widehat{Q}_n=\widehat{U}\biggl(\frac{n}{k}\biggr)\biggl(\frac{k\widehat{\nu
(S)}}{np}\biggr)^{{1}/{\widehat{\alpha
}}}\widehat{S}
\end{eqnarray}
with
\[
\widehat{S}= \bigl\{\mathbf{z}=r\mathbf{w}\dvtx  r\geq(\widehat{\alpha
}\widehat{\psi
}_n(\mathbf{w}))^{{1}/({\widehat{\alpha }+d})}, \mathbf{w}\in
\Theta\bigr\}
\]
and
\[
\widehat{\nu(S)}=\widehat{\alpha }^{-{\widehat{\alpha
}}/{(\widehat{\alpha }+d})}\int_{\Theta}
(\widehat{\psi}_n(\mathbf{w}))^{{d}/{(\widehat{\alpha
}+d)}}\,d\lambda(\mathbf{w}).
\]
In the definition of the set $S$, the choice of the value 1 was not
motivated. We could have taken any number $c>0$ instead.
Such an alternative definition of $S$ would lead to exactly the same
estimator $\widehat{Q}_n$, which shows that the value 1 plays no role.

Assume
{\renewcommand{\theequation}{c}
\begin{equation}\label{eqc}
\lim_{t\rightarrow\infty }\frac{U(t)}{t^{1/\alpha}}=c\qquad  \mbox{for
some } c\in(0,\infty).
\end{equation}}
\vspace*{-\baselineskip}\vspace*{10pt}
\setcounter{equation}{9}

\noindent Note that this simple condition is weaker than the usual second order
condition with negative second order parameter $\rho$ [see, e.g.,
Theorem 4.3.8 in de Haan and Ferreira (\citeyear{deHFer06})]; indeed,
there exist functions $U$ with $\rho=0$ that satisfy condition~(\ref{eqc}).

\begin{theorem}\label{thm1}
Let $p\to0$ as $n\rightarrow\infty $.
Assume conditions \textup{(a), (b), (\ref{eqc})} hold and that $\widehat{\alpha }$ is
such that
$\sqrt{k}(\widehat{\alpha }-\alpha
)=O_{\mathbb{P}}(1)$.
Also assume that $({\log np})/{\sqrt{k}}\rightarrow0$, \mbox{$h\to0$} and
$k/(c(h, K)\log k)\to\infty$, as $n\to\infty$.
Then we have
\begin{equation} \label{(theorem)}
\frac{P(\widehat{Q}_n\triangle Q_n)}{p} \stackrel{\mathbb{P}}{\rightarrow}0\qquad
\mbox{as }
n\rightarrow\infty ,
\end{equation}
and hence
\[
\frac{P(\widehat Q_n)}{ p}\stackrel{\mathbb{P}}{\to} 1.
\]
\end{theorem}

\begin{remark}\label{rem1}
The tuning parameter $k$ is used in the estimators of $\alpha,
U(n/{k})$ and $\psi$. It is important to be able to choose three
different values for $k$, denoted with $k_\alpha , k_U$ and $k_\psi$,
respectively. (Note that ``good'' values of $k_\alpha $ and $k_U$ are
determined by the tail of $H$---the distribution function of $R_1$---whereas a good $k_\psi$ is determined by the conditional distribution
of $\mathbf{W}_1$, given that $R_1>r$, for large $r$.) If we adapt the
conditions of the theorem, in particular if (b) holds for $k_\alpha , k_U,
k_\psi$ and if $(\log np) / \sqrt{k_\alpha } \to0$, $k_\psi/(c(h,
K)\log
k_\psi)\to\infty$ and $(\log k_U)/\sqrt{k_\alpha }\to0$, then~(\ref{(theorem)}) remains true for the generalized estimator that
allows for the aforementioned different $k$-values. We will use this
generalized estimator in the simulation study and the real data
application.

The actual choice of these $k$-values is a notorious problem in extreme
value theory. A solution of this problem is far beyond the scope of the
present paper. We will only give heuristic guidelines here. First,
consider the estimation of $\alpha $. Plot $\widehat{\alpha }$ as a
function of $k$. Now
find the first stable, that is,\ approximately constant, region in the
graph of this function. This vertical level is the final estimate of
$\alpha $. It is also possible to use (complicated) \textit{asymptotically}
optimal procedures; see, for example, Danielsson et al.\ (\citeyear
{Danetal01}). Once the
estimate $\widehat{\alpha }$ is fixed, we plot $\widehat{U}(\frac
{n}{k})(\frac{k}{n})^{1/\widehat{\alpha }}$
against $k$ and we search for the first stable part in this graph. The
vertical level is now the estimate of the constant $c$ in condition
(c). Observe that $\widehat{U}(\frac{n}{k})(\frac{k}{n})^{1/\widehat
{\alpha }}$ is a building block of
$\widehat Q_n$, so we do not need to estimate $U(\frac{n}{k})$
separately. Also
observe that we do not (need to) determine $k_\alpha $ and $k_U$, but
only a
region of good values. Finally, using again the already fixed $\widehat
{\alpha }$, we
plot $\widehat{\nu(S)}$ as a function of $k$ and again we search for
the first stable region; we take $k_\psi$ to be the midpoint of this
region of $k$-values.
\end{remark}

\begin{remark}\label{rem2}
The class of multivariate regularly varying distributions is quite
large. It contains, for example, all elliptical distributions with a heavy
tailed radial distribution and all distributions in the domain of a
sum-attraction of a multivariate (nonnormal) stable distribution. It
seems natural, however, to try to extend the assumption of multivariate
regular variation to the case of nonequal tail indices~$\alpha$. It is
an important feature of the present model that all directions are
equally important: the marginal distributions do not play a special
role. An extension to nonequal tail indices would be possible in
principle, but it will be of limited value since it only works if
\textit{marginal} transformations lead to the present model. Also note
that basically all linear combinations of the components inherit the
lowest of the marginal tail indices: the tail index is not a smooth
function of the direction (if it is not constant). Moreover, the
statistical theory that will be needed will be challenging and will
lead to a new and different project.
\end{remark}

\begin{remark}\label{rem3}
Note that the estimated extreme risk region $\widehat Q_n=\widehat{Q}_n(p)$
depends on $p$ in a continuous way and has the property that $p_1<p_2$
implies $\widehat{Q}_n(p_1)\subset\widehat{Q}_n(p_2)$. Hence, we can
find the smallest $p$
such that an observation is on the boundary of $\widehat Q_n(p)$. The
corresponding observation can be considered the largest one and we
know its ``$p$-value.'' This is helpful in deciding whether some
observation is the most extreme or if it is an outlier. Also, by
continuing this procedure we can rank the larger observations.
\end{remark}

\section{Simulation study}\label{sec3}
In this section, a detailed simulation study is performed in order to
investigate the finite sample performance of our estimator [with
$1/\alpha$ estimated using the moment estimator of \citet
{DekEindeH89} and with $K(u)=1-u$].
We consider five multivariate distributions.
\begin{itemize}

\item The bivariate Cauchy distribution with
density
\begin{equation}\label{cauchy2}f(x,y)=\frac{1}{2\pi
(1+x^2+y^2)^{3/2}} ,\qquad  (x,y) \in\mathbb{R}^2.
\end{equation}
This is a very heavy tailed density, with $\alpha=1$ and
$\psi(\mathbf{w})=1/(2\pi)$, for $\mathbf{w}\in\Theta$.
\item The trivariate Cauchy distribution with
density
\begin{equation}\label{cauchy3}f(x,y,z)=\frac{1}{\pi
^2(1+x^2+y^2+z^2)^2} ,\qquad  (x,y, z) \in\mathbb{R}^3 .\vadjust{\goodbreak}
\end{equation}
This is also a very heavy tailed density, with $\alpha=1$
and $\psi(\mathbf{w})=1/(4\pi)$, for $\mathbf{w}\in\Theta$.
\item A bivariate elliptical distribution with density ($r_0\approx1.2481$)
\begin{equation}\label{ellip}
f(x,y)= \cases{
\displaystyle \frac{3}{4\pi}r_0^4(1+r_0^6)^{-3/2}, &\quad $x^2/4+y^2<
r_0^2,$\vspace*{2pt}\cr
\displaystyle \frac{3(x^2/4+y^2)^2}{4\pi(1+(x^2/4+y^2)^3)^{3/2}} ,&\quad $ x^2/4+y^2\geq
r_0^2. $
}
\end{equation}
%
It is less heavy tailed. We have $\alpha=3$ and
$\psi(w_1,w_2)=c(1+3w_2^2)^{-5/2}, \mathbf{w}=(w_1,w_2)\in\Theta$, with
$c\approx0.6028$.
\item A bivariate ``clover'' distribution with density ($r_0\approx1.2481$)
\begin{equation}\label{clov}
f(x,y)= \cases{
\displaystyle \frac{3}{10\pi}r_0^4(1+r_0^6)^{-3/2} \biggl(5+\frac
{4(x^2+y^2)^2-
32x^2y^2}{r_0(x^2+y^2)^{3/2}} \biggr),\vspace*{3pt}\cr
\qquad   x^2+y^2< r_0^2,\vspace*{4pt}\cr
\displaystyle \frac{3 (9(x^2+y^2)^2-32x^2y^2 )}{10\pi(1+({x^2}+y^2)^3)^{3/2}}
,\vspace*{3pt}\cr
\qquad x^2+y^2\geq r_0^2.
}\hspace*{-8pt}
\end{equation}
%
We have $\alpha=3$, again, and
$\psi(w_1,w_2)=(9-32w_1^2w_2^2)/(10\pi), \mathbf{w}=(w_1,\break w_2)\in
\Theta$.
\item A bivariate asymmetric shifted distribution with density
[$r_0\approx1.2331$, $\widetilde{r}(x,y):=r_0\vee((x+5)^2+y^2)^{1/2}$]
\begin{equation}\label{asymm}
f(x,y) = \cases{
\displaystyle \frac{\widetilde{r}^2(x,y)}{6\pi(1+\widetilde
{r}^4(x,y))^{5/4}} \biggl(3+\frac{x+5}{\widetilde{r}(x,y)} \biggr),
\hspace*{69pt}\qquad y\geq0, \vspace*{5pt}\cr
\displaystyle \frac{\widetilde{r}^2(x,y)}{6\pi(1+\widetilde{r}^4(x,y))^{5/4}}
\biggl(3+\frac{(x+5)^3-3(x+5)y^2}{\widetilde{r}^3(x,y)} \biggr),
\qquad   y<0.
}\hspace*{-15pt}
\end{equation}
%
This distribution is not symmetric and the ``center'' is not the origin,
but $(-5,0)$; $\alpha=1$ and
$\psi(w_1,w_2)= \frac{1}{6\pi}(3+w_1)$, if $ w_2\geq0$, and
$ \psi(w_1,w_2)= \frac{1}{6\pi}(3+4w_1^3-3w_1)$, if $w_2<0$,
$\mathbf{w}=(w_1,w_2)\in\Theta$.
\end{itemize}

\begin{figure}

\includegraphics{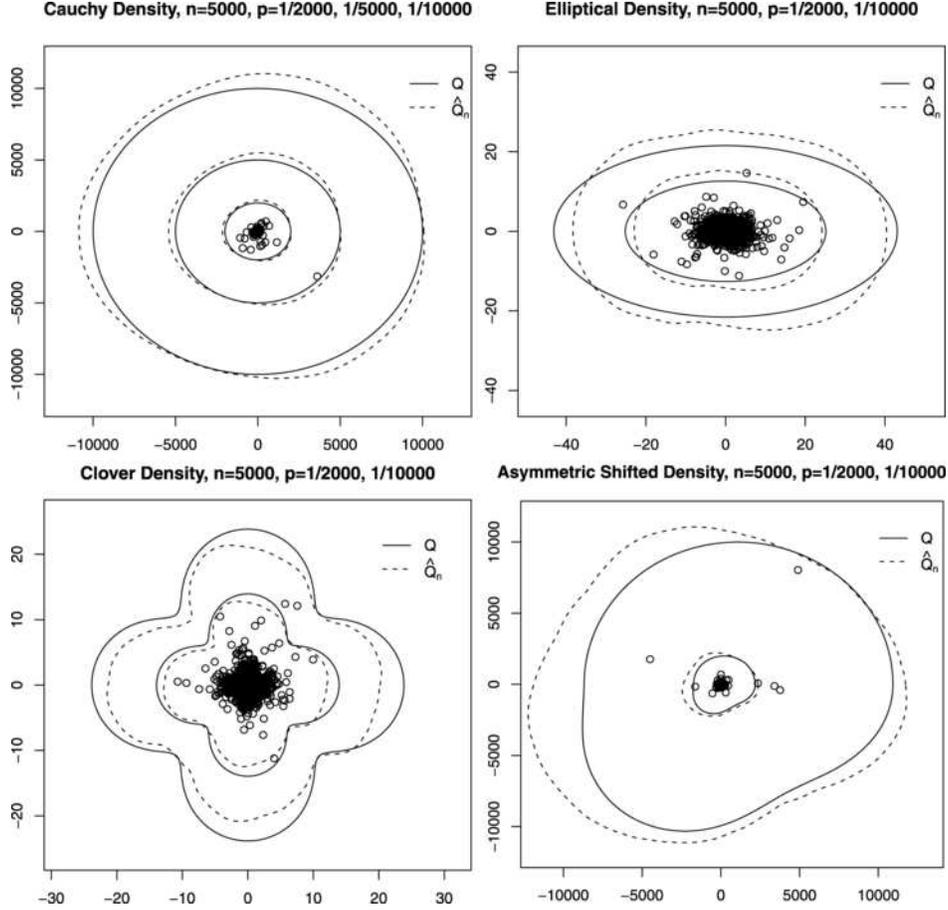}

\caption{True and estimated risk regions based on
one sample of size 5,000 from the bivariate Cauchy distribution, the
elliptical distribution in (\protect\ref{ellip}), the clover distribution in
(\protect\ref{clov}) and the asymmetric shifted distribution in (\protect\ref{asymm}).}\label{clouds}
\end{figure}

First, we simulated single data sets of size 5,000 of the bivariate
Cauchy distribution, the elliptical distribution in (\ref{ellip}), the
clover distribution in (\ref{clov}) and the asymmetric shifted
distribution in (\ref{asymm}). We computed the true and estimated risk regions
for $p=1/2\mbox{,}000$, 1$/$5,000 or $1/10\mbox{,}000$. This is depicted in Figure \ref
{clouds}. We see that the estimated regions are relatively close to the
true risk regions. It is interesting to note that the $p$-value (see
Remark~\ref{rem3}) of the largest observation for the Cauchy sample is
0.000209, which is about $1/n$. This shows that this observation is a
typical one. (Looking at the data only, one might want to conclude that
this observation is an outlier.) Also note that for the bivariate
Cauchy distribution, for, for example, $p=1/10\mbox{,}000$, the density $f$ at the
boundary of the true risk region is less than $10^{-12}$. This
emphasizes that we are estimating in an ``almost empty'' part of the
plane and that a fully nonparametric procedure could not work here.

\begin{figure}

\includegraphics{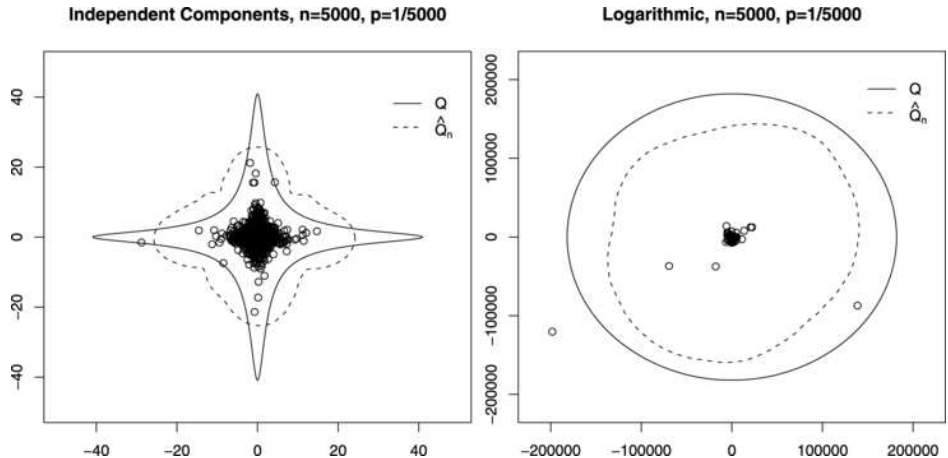}

\caption{True and estimated risk regions based on
one sample of size $5\mbox{,}000$ from the bivariate distribution with
independent $t_3$-components and the ``logarithmic'' distribution.}\label{clouds2}
\end{figure}

In addition, we simulated one sample of the bivariate distribution with
independent $t_3$-components. This distribution does \textit{not}
satisfy condition (a), since the spectral measure is discrete and
concentrated on the intersection of the coordinate axes with the unit
circle. We also simulated one sample of a bivariate ``logarithmic''
distribution with $\alpha=1$ and uniform spectral measure, but where
the radial distribution satisfies $U(t)/(t\log t)$ tends to a constant
and hence $U(t)/t \to\infty$ as $t\to\infty$, that is,\ this
distribution does \textit{not} satisfy condition (c). Although both
distributions do not satisfy our conditions, we see nevertheless
satisfactory behavior of the estimator in Figure~\ref{clouds2}. In the
left panel, the estimated region has about the right size and the
difficult shape is approximated reasonably well; in the right panel, we
see that both the shape and the size are approximated quite well.

%

After this visual assessment of our estimator based on one sample at a
time, we now investigate its performance based on 100 simulated samples
of size 5,000. We will compare our estimator (denoted EVT) to a
nonparametric and to a more parametric estimator. The nonparametric
estimator is only defined in case $p=1/n$ and tries to mimic the
largest order statistic as an estimator of the $(1-1/n)$th quantile in
the univariate case. It aims at elliptical level sets. It is defined as
follows. First, calculate the smallest ellipsoid containing half of the
data, the so-called MVE. Then inflate this ellipsoid, such that the
``largest'' observation lies on its boundary. Now the region
outside this ellipsoid is the estimator.

For $d=2$, the more parametric estimator is defined similarly to
$\widehat{Q}_n$
in (\ref{defqhat}), but (only) the estimation of $(\nu(S))^{1/\alpha
}S$ is done parametrically. Therefore, this estimator has the same size
as $\widehat{Q}_n$, but a different shape. (Note that the fully parametric
estimator based on multivariate normality would have a very bad
performance.) Take the $k$ observations with radius $R_i> R_{n-k:  n}$
and consider the transformed data
$(R_i/R_{n-k:  n}, \mathbf{W}_i)$. 
In line with the limit result in (\ref{mrv}), assume that these
data have a ``distribution'' $(\cdot)^{-\alpha}\Psi$, where $\Psi$
depends on a parameter $\rho$. To be precise, we assume for the density
\[
\psi_\rho(\theta)=(4\pi)^{-1}\bigl(2+\sin\bigl(2(\theta-\rho)\bigr)\bigr),\qquad  0\leq
\theta< 2\pi, 0\leq\rho<\pi.
\]
(Here a point on the unit circle is represented by its angle $\theta
\in[0,2\pi)$.) Now $\alpha$ and $\rho$ are estimated by maximum
likelihood; observe that this yields the Hill estimator for $1/\alpha$.

\begin{table}
\tablewidth=265pt
\caption{Median of the relative errors $P(\widehat
{Q}_n\triangle Q_n)/p$
of the three estimators, for $p=1/5\mbox{,}000$ (p1) and $1/10\mbox{,}000$ (p2)}
\label{tab}
\begin{tabular*}{\tablewidth}{@{\extracolsep{\fill}}lccccc@{}}
\hline
\textbf{Density}& \textbf{EVT p1} & \textbf{Par p1} & \textbf{NP p1} & \textbf{EVT p2} & \textbf{Par p2} \\
\hline
Biv.\ Cauchy & $0.28$ & $0.29$ & $0.72$ & $0.31$ & $0.32$ \\
Triv.\ Cauchy & $0.22$ & -- & $0.54$ & $0.24$ & -- \\
Elliptical & $0.36$ & $0.51$ & $0.80$ & $0.39$ & $0.54$ \\
Clover & $0.44$ & $0.57$ & $0.58$ & $0.49$ & $0.61$\\
Asymm.\ shifted & $0.26$ & $0.27$ & $0.61$ & $0.30$ & $0.32$\\
\hline
\end{tabular*}       \vspace*{-3pt}
\end{table}

\begin{figure}

\includegraphics{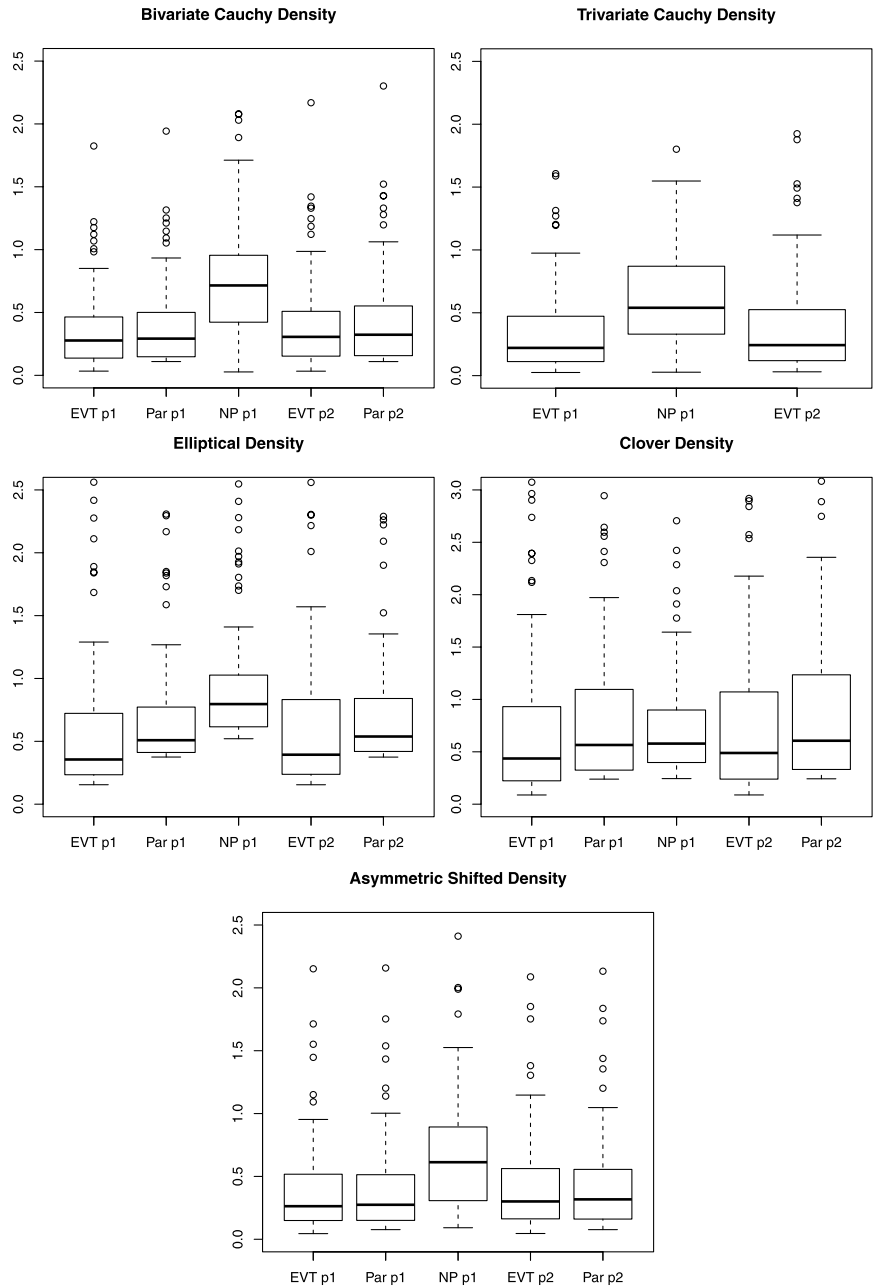}

\caption{Boxplots of ${P(\widehat{Q}_n\triangle Q_n )}/{p}$
for the here
proposed estimator and for the parametric and the nonparametric
estimator, based on
$100$ simulated data sets of size 5,000 from the five presented
densities for $p=1/5\mbox{,}000$ (p1) and $1/10\mbox{,}000$ (p2).}\label{f3}
\end{figure}

Table~\ref{tab} shows for the three different estimators the median
of the 100 relative errors $P(\widehat{Q}_n\triangle Q_n)/p$ for
$p=1/5\mbox{,}000$ (p1)
and $1/10\mbox{,}000$ (p2).
In Figure~\ref{f3}, boxplots are shown of the relative error
$P(\widehat{Q}_n\triangle Q_n)/p$ for $p=1/5\mbox{,}000$ (p1) and $1/10\mbox{,}000$
(p2). From
this table and figure, we see a good performance of our estimator. Its
behavior does not change much if $p$ changes from $1/5\mbox{,}000$ to $1/10\mbox{,}000$.
The parametric estimator performs reasonably well, but it is
outperformed by our estimator, in particular for the elliptical and
clover densities. Recall that this estimator can be seen as a
modification of our estimator, since it uses the same estimated
inflation factor, but the shape is estimated differently. We see a
moderate performance of the nonparametric estimator; also, it cannot be
adapted to $p=1/10\mbox{,}000$. Given that the estimation of these extreme
risk regions is a statistically difficult problem, we see decent
behavior of the three estimation methods. Obviously the parametric and
the nonparametric estimator do not perform well if the parametric part
of the model is not adequate or if the shape of the region is not
elliptical, respectively. The EVT estimator, presented in this paper,
does not suffer from these shortcomings and performs well for many
multivariate distributions.

\begin{figure}

\includegraphics{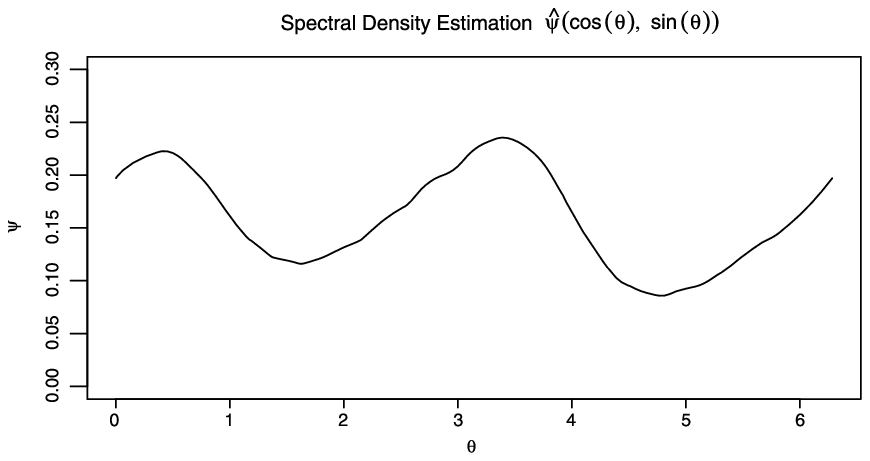}

\caption{Estimator of $\psi$ of bivariate exchange rate
returns.}\label{psi}\vspace*{-3pt}
\end{figure}

\begin{figure}[b]

\includegraphics{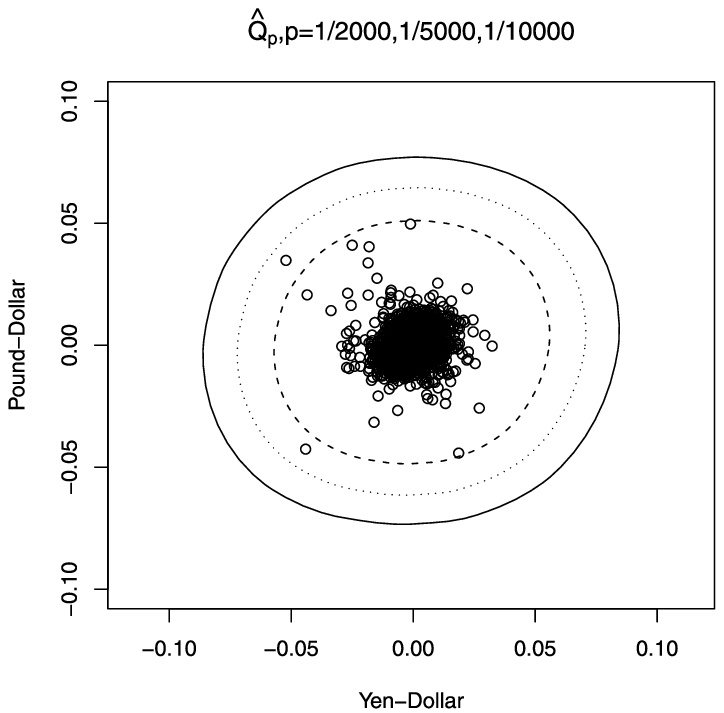}

\caption{Estimated extreme risk regions of exchange rate
returns.}\label{66}
\end{figure}

\section{Application}\label{sec4}
In this section, an application of our method to foreign exchange rate
data is presented. The data
are the daily exchange rates of Yen-Dollar and Pound-Dollar from
January 4, 1999 to July 31, 2009. Consider the daily log returns given
by
$
X_{i,j}=\log{({Y_{i+1,j}}/{Y_{i,j}})}$,
with $i=1,\ldots,2\mbox{,}664$, $j=1,2$, and $Y_{i,1}$ is the daily
exchange rate of the Yen to the Dollar and $Y_{i,2}$ of the Pound to
the Dollar.
First, we check the equality of the extreme value indices (the
reciprocals of the tail indices) of the right and left tails of both
marginal distributions and that of the radius. This yields 5 extreme
value indices; the 5 estimates in increasing order are: 0.141, 0.191,
0.223, 0.242, 0.256. Hence, the maximal difference is 0.115. Based on
the asymptotic normality of\vadjust{\goodbreak} the moment estimator of the extreme value
index, we compute an approximate upper bound for the maximal difference
of the 5 estimators under the null hypothesis of equality: 0.264.
Hence, there is no evidence that the 5 extreme value indices are
different. Other exchange rate data sets share this property. There are
also economic arguments supporting this claim. Therefore, we estimate
$\alpha$ based on the radius and find $\widehat\alpha=3.90$.
As a next step, we estimate the density $\psi$ of the spectral
measure. The estimate $\widehat\psi_n$
is depicted in Figure 4;
it is almost periodic with period $\pi$. This yields that the boundary
of the
estimated extreme risk region is not like a circle, but more like an
ellipse. The location of the maxima of $\widehat\psi_n$ correspond to
the major axis of the region. We estimate the extreme risk regions for
$p=1/2\mbox{,}000,1/5\mbox{,}000$ and $1/10$,$000$; see Figure\vadjust{\goodbreak} 5. 
For risk management of financial institutions in the U.S.,\ it is
important to know which extreme exchange rate returns w.r.t.\ the Pound
and the Yen can occur and which returns essentially never occur. Our
estimate answers this question.
More specifically,
points on the boundary of the estimated extreme risk region can be used
as multivariate stress test scenarios.
A scenario on the intersection of the major axis of the
ellipse-like boundary of the extreme risk region and the boundary
itself corresponds to a larger shock than a scenario on the
intersection of the minor axis of the ellipse-like boundary and the
boundary itself, but our method shows
that their ``extremeness'' is about the same.
%

\section{Proofs}\label{sec5}

For the proof of the theorem, we need several lemmas and propositions.
We assume throughout that the conditions of the theorem are in force.
We start with a lemma on regular
variation in $\mathbb{R}^d$.

\begin{lemma}\label{lem1}
Write $l=1/\int_{\{\|\mathbf{z}\|\geq1\}
}q(\mathbf{z}
)\,d\mathbf{z}$. For any $\varepsilon >0$,
\begin{equation}\label{(4.1)}
\lim_{t\rightarrow\infty }\sup_{\|\mathbf{z}\|\geq\varepsilon }
\biggl|\frac{f(t\mathbf{z}
)}{t^{-d}V(t)}-q(\mathbf{z}) \biggr|=0.
\end{equation}
Moreover
\begin{equation}\label{ne}
\lim_{t\rightarrow\infty }\frac{\mathbb{P}(\mathbf{X}\in
tB)}{V(t)}=\int
_{B}q(\mathbf{z})\,d\mathbf{z}
\end{equation}
for any Borel set $B$ bounded away from the origin.
Define $q_t(\mathbf{z})=t(U(t))^d\times f(U(t)\mathbf{z})$. Then
\begin{eqnarray}\label{(4.4)}
\lim_{t\rightarrow\infty }\sup_{\|\mathbf{z}\|\geq\varepsilon }
|q_t(\mathbf{z})-lq(\mathbf{z})|=0.
\end{eqnarray}
Let $\widetilde h$ be the density of $H$, then
\begin{equation}\label{(4.5)}
\lim_{t\rightarrow\infty }\frac{\widetilde h(t)}{t^{-1}V(t)}=\frac
{\alpha}{l} .
\end{equation}
\end{lemma}

\begin{pf}
For any $\|\mathbf{z}\|\geq\varepsilon >0$
[cf. Theorem 2.1 in de Haan and Resnick (\citeyear{deHRes87})],
\begin{eqnarray*}
&& \biggl|\frac{f(t\mathbf{z})}{t^{-d}V(t)}-q(\mathbf{z}) \biggr|\\[-2pt]
&&\qquad = \biggl|\frac{f(t\|\mathbf{z}\|(\mathbf{z}/{\|\mathbf{z}
\|}))}{(t\|z\|)^{-d}V(t\|z\|)}\cdot\frac{(t\|z\|)^{-d}V(t\|z\|)}{t^{-d}V(t)}-
q(\mathbf{z}) \biggr|\\[-2pt]
&&\qquad \leq\|\mathbf{z}\|^{-d-\alpha } \biggl|\frac{f(t\|\mathbf{z}\|
({\mathbf{z}}/{\|\mathbf{z}
\|}))}{(t\|z\|)^{-d}V(t\|z\|)}
-q \biggl(\frac{\mathbf{z}}{\|\mathbf{z}\|} \biggr) \biggr|\\[-2pt]
&&\qquad \quad {}
+\frac{f(t\|\mathbf{z}\|({\mathbf{z}}/{\|\mathbf
{z}\|}))}{(t\|z\|)^{-d}V(t\|z\|)}
\biggl|\frac{(t\|z\|)^{-d}V(t\|z\|)}{t^{-d}V(t)}-\|\mathbf{z}\|^{-d-\alpha
} \biggr|
\end{eqnarray*}
Then ({\ref{(4.1)}}) follows from condition (a).\vadjust{\goodbreak}

Let a Borel set $B$ be such that
$B\subset\{\|\mathbf{z}\|\geq\gamma\}$, for some $\gamma>0$. Then
for $\mathbf{z}\in B$ and sufficiently large $t$,
${f(t\mathbf{z})}/{t^{-d}V(t)}$ is bounded by $q(\|\mathbf
{z}\|^{-1}\mathbf{z})\|\mathbf{z}\|^{-a/2-d}$.
Hence, (\ref{ne}) holds by Lebesgue's dominated convergence
theorem.

We have from (\ref{ne}), as $t\rightarrow\infty $,
\[
tV(U(t))=\frac{V(U(t))}{\mathbb{P}(R\geq U(t))}\rightarrow l.
\]
Hence (\ref{(4.1)}) implies, uniformly for $\| \mathbf{z}\|
\geq\varepsilon$,
\[
q_t(\mathbf{z})=tV(U(t))\frac{f(U(t)\mathbf
{z})}{(U(t))^{-d}V(U(t))}\rightarrow l
q(\mathbf{z}).
\]

Note that
\[
1-H(t)=\mathbb{P}(R> t)=\int_t^\infty \int_{\Theta}f(r\mathbf{w})\,d\lambda
(\mathbf{w})r^{d-1}\,dr.
\]
By taking derivatives, (\ref{(4.1)}) and the homogeneity of $q$, we obtain
\[
\lim_{t\rightarrow\infty }\frac{\widetilde h(t)}{t^{-1}V(t)}
=\int_{\Theta}q(\mathbf{w})\,d\lambda (\mathbf{w})
=\alpha \int_{\{\|\mathbf{z}\|\geq1\}}q(\mathbf{z})\,d\mathbf
{z}=\alpha /l.
\]
\upqed
\end{pf}

We now see that (\ref{pw}) and (\ref{(1.3)}) hold with $V=1-H$. From
now on, we will make the choice $V=1-H$ and hence $l=1$.
Note that with this choice the relations (\ref{(1.1)}) [with $\nu
(B)=\int_B q(\mathbf{z})\,d\mathbf{z}$] and (\ref{(1.2)}) readily
follow from ({\ref{ne}}).

\begin{coro}\label{coro1}
For all Borel sets $B$ with
positive distance from the origin,
\begin{equation}\label{(4.3)}
\lim_{t\rightarrow\infty }tP(U(t)B)=\nu(B)
\end{equation}
and
\begin{equation}\label{(4.32)}
\lim_{n\rightarrow\infty } \frac{\nu(S)}{p}P\biggl(U \biggl( \frac{n}{k}
\biggr)\biggl(\frac
{k\nu(S)}{np}\biggr)^{1/{\alpha }}B\biggr)=\nu(B).
\end{equation}
\end{coro}

\begin{pf}
From $\mathbb{P}(R\geq U(t))=1/t$ and (\ref{(1.1)}), we
obtain (\ref{(4.3)}). 
It
follows from (c)
that
\begin{equation}\label{(uvsp)}
\frac{U({\nu(S)}/{p})}{U({n}/{k})({k\nu
(S)}/{(np)})^{{1}/{\alpha }}}
\rightarrow1.
\end{equation}
This yields (\ref{(4.32)}).
\end{pf}

\begin{lemma}\label{lem2}
For each $\varepsilon >0$, there exists a
$\delta >0$ and $t_0>0$ such that for $t>t_0$
\[
\biggl\{\mathbf{z}\dvtx \frac{f(t\mathbf{z})}{t^{-d}V(t)}\leq\varepsilon \biggr\}
\subset
\{\mathbf{z}\dvtx\|\mathbf{z}\|>\delta \}.
\]
\end{lemma}

\begin{pf}
It is sufficient to prove
$\{\mathbf{z}\dvtx \|\mathbf{z}\|\leq\delta \}\subset\{\mathbf
{z}\dvtx f(t\mathbf{z})/(t^{-d}V(t))>\varepsilon \}$.
First, by (\ref{(1.3)}) and the continuity of $q$, for some $c_1>0$,
there exists $s_0>0$ such that
for $s>s_0$
\[
\inf_{\mathbf{w}\in\Theta}\frac{f(s\mathbf{w})}{s^{-d}V(s)}\geq c_1
\]
and also for $s_1, s_2>s_0$ [cf.\ Proposition B.1.9.5 in de
Haan and Ferreira (\citeyear{deHFer06})]
\[
\frac{V(s_1)}{V(s_2)}>\frac{1}{2}\biggl(\frac{s_1}{s_2}\biggr)^{-\alpha /2}.
\]
Now for $t>s_0$ and any $\mathbf{z}\in\{\mathbf{z}\dvtx \|\mathbf
{z}\|\leq\delta \}$, there
are two possibilities.
\begin{longlist}[(ii)]
\item[(i)] $t\|\mathbf{z}\| > s_0$, then
\[
\frac{f(t\mathbf{z})}{t^{-d}V(t)}
=\frac{f(t\|\mathbf{z}\|(\mathbf{z}/{\|\mathbf
{z}\|}))}{(t\|\mathbf{z}\|)^{-d}V(t\|\mathbf{z}
\|)}\cdot\frac{(t\|\mathbf{z}\|)^{-d}V(t\|\mathbf{z}\|)}{t^{-d}V(t)}
> \frac{1}{2}c_1 \delta ^{-\alpha /2-d}>\varepsilon ;
\]
\item[(ii)] $t\|\mathbf{z}\| \leq s_0$, then by continuity of
$f$ and $f>0$, we have for some $c_2>0$, $f(t\mathbf{z})\geq
c_2$, and hence, since $\lim_{t\to\infty }t^{-d}V(t)=0$, we
obtain for $t>t_0(\geq s_0)$
\[
\frac{f(t\mathbf{z})}{t^{-d}V(t)}>\varepsilon .
\]
\end{longlist}
\upqed
\end{pf}

\begin{lemma}\label{lem3}
For $\varepsilon >0$ and large $n$,
\[
\bar{Q}_n\subset U \biggl(\frac{\nu(S)}{p} \biggr)\{\mathbf{z}\dvtx  q(\mathbf
{z})\leq1+\varepsilon \}
\]
and
\[
\bar{Q}_n\supset U \biggl(\frac{\nu(S)}{p} \biggr)\{\mathbf{z}\dvtx  q(\mathbf
{z})\leq1-\varepsilon \}.
\]
\end{lemma}

\begin{pf}
Recall that
$\bar{Q}_n= \{\mathbf{z}\dvtx f(\mathbf{z})\leq(\frac{np}{k\nu
(S)})^{1+({d}/{\alpha})}\frac{1}{({n}/{k})(U({n}/{k}))^d} \}$.
It follows from (\ref{(uvsp)}) that for $n$ large enough and
$\varepsilon _1>0$
\begin{eqnarray*}
\bar{Q}_n&=&U\biggl(\frac{\nu(S)}{p}\biggr)
\biggl\{\mathbf{z}\dvtx  f\biggl(U\biggl(\frac{\nu(S)}{p}\biggr)\mathbf{z}\biggr)\leq\biggl(\frac{np}{k\nu
(S)}\biggr)^{1+
({d}/{\alpha })}\frac{1}{({n}/{k})(U({n}/{k}))^d} \biggr\}\\
&=&U\biggl(\frac{\nu(S)}{p}\biggr)
\biggl\{\mathbf{z}\dvtx  q_{\nu(S)/p}(\mathbf{z})\leq\biggl(\frac{np}{k\nu
(S)}\biggr)^{{d}/{\alpha }}
\biggl(U\biggl(\frac{n}{k}\biggr)\biggr)^{-d}\biggl(U\biggl(\frac{\nu(S)}{p}\biggr)\biggr)^{d} \biggr\}\\
&\subset& U\biggl(\frac{\nu(S)}{p}\biggr) \bigl\{\mathbf{z}\dvtx  q_{\nu(S)/p}(\mathbf
{z})\leq1+\varepsilon _1\bigr\}.
\end{eqnarray*}
Now Lemma~\ref{lem2} implies
$ \{\mathbf{z}\dvtx  q_{\nu(S)/p}(\mathbf{z})\leq1+\varepsilon_1 \}
\subset\{\mathbf{z}\dvtx \|\mathbf{z}
\|>\delta \}$,
hence we have by (\ref{(4.4)})
\[
\bar{Q}_n\subset U \biggl(\frac{\nu(S)}{p} \biggr)\{\mathbf{z}\dvtx  q(\mathbf
{z})\leq1+\varepsilon \}.
\]
The other inclusion follows in the same way (but Lemma~\ref{lem2} is not
needed).
\end{pf}

\begin{lemma}\label{lem4}
For $\varepsilon >0$ and large $n$,
\[
\widetilde{Q}_n\subset U \biggl(\frac{\nu(S)}{p} \biggr)\{\mathbf{z}\dvtx  q(\mathbf
{z})\leq1+\varepsilon \}
\]
and
\[
\widetilde{Q}_n\supset U \biggl(\frac{\nu(S)}{p} \biggr)\{\mathbf{z}\dvtx  q(\mathbf
{z})\leq1-\varepsilon \}.
\]
\end{lemma}

\begin{pf}
Recall that
$\widetilde{Q}_n=U ( \frac{n}{k} )(\frac{k\nu(S)}{np})^{
{1}/{\alpha }}\{\mathbf{z}\dvtx q(z)\leq1\}$.

Put $T_n=(U(\frac{\nu(S)}{p}))^{-1}U (\frac{n}{k} )(\frac{k\nu
(S)}{np})^{{1}/{\alpha }}$,
then
\begin{eqnarray*}
\widetilde{Q}_n&=&U\biggl(\frac{\nu(S)}{p}\biggr) \{T_n\mathbf{z}\dvtx q(\mathbf
{z})\leq1 \}
=U\biggl(\frac{\nu(S)}{p}\biggr) \{T_n\mathbf{z}\dvtx q(T_n\mathbf{z})\leq
T_n^{-d-\alpha } \}\\
&=&U\biggl(\frac{\nu(S)}{p}\biggr) \{\mathbf{z}\dvtx q(\mathbf{z})\leq
T_n^{-d-\alpha } \}.
\end{eqnarray*}
Since $T_n\rightarrow1$ as $n\rightarrow\infty $ by
(\ref{(uvsp)}), the result follows.
\end{pf}

\begin{prop}\label{prop1}
We have
\[
\lim_{n\rightarrow\infty }\frac{P(Q_n\bigtriangleup\widetilde{Q}_n)}{p}=0.
\]
\end{prop}

\begin{pf}
Note that $P(Q_n\bigtriangleup\widetilde
{Q}_n)\leq
P(Q_n\bigtriangleup\bar{Q}_n)+P(\bar{Q}_n\bigtriangleup\widetilde
{Q}_n)$. Observe that
$Q_n\subset\bar{Q}_n$ or $\bar{Q}_n\subset Q_n$, hence
$P(Q_n\bigtriangleup\bar{Q}_n
)\leq|p-P(\bar{Q}_n)|$. By Lemma~\ref{lem3} and Corollary~\ref{coro1},
for any $\varepsilon >0$ and large $n$
\begin{eqnarray*}
\frac{\nu(S)}{p}P(\bar{Q}_n)&\leq& \frac{\nu(S)}{p}P\biggl(U\biggl(\frac{\nu
(S)}{p}\biggr)\{\mathbf{z}\dvtx q(\mathbf{z})\leq1+\varepsilon \}\biggr)\\
&\rightarrow& \nu\bigl(\{\mathbf{z}\dvtx q(\mathbf{z})\leq1+\varepsilon \}
\bigr)\\
&=&\nu\bigl(\bigl\{\mathbf{z}\dvtx  q\bigl(\mathbf{z}(1+\varepsilon )^{1/(d+\alpha
)}\bigr)\leq1\bigr\}\bigr)\\
&=&\nu\bigl(\bigl\{(1+\varepsilon )^{-1/(d+\alpha )}\mathbf{z}\dvtx q(\mathbf
{z})\leq1\bigr\}\bigr)\\
&=&(1+\varepsilon )^{\alpha /(d+\alpha )}\nu(S).
\end{eqnarray*}
Thus, $\limsup_{n\rightarrow\infty }\frac{P(\bar{Q}_n)}{p}\leq
(1+\varepsilon )^{\alpha
/(2+\alpha )}$.\vadjust{\goodbreak}

Similarly, we have $\liminf_{n\rightarrow\infty }\frac{P(\bar
{Q}_n)}{p}\geq
(1-\varepsilon )^{\alpha /(2+\alpha )}$.
Hence,\break $\lim_{n\rightarrow\infty }\frac{P(\bar{Q}_n)}{p}=1$, that
is,\
$\lim_{n\rightarrow\infty }\frac{P(Q_n\bigtriangleup\bar{Q}_n)}{p}=0$.

In the same way, it follows from Lemmas~\ref{lem3} and~\ref{lem4} that
\begin{eqnarray*}
\frac{\nu(S)}{p}P(\bar{Q}_n\bigtriangleup\widetilde{Q}_n)
&\leq& \frac{\nu(S)}{p}P\biggl(U\biggl(\frac{\nu(S)}{p}\biggr)\{\mathbf
{z}\dvtx 1-\varepsilon \leq q(\mathbf{z}
)\leq1+\varepsilon \}\biggr)\\
&\rightarrow& \nu\bigl(\{\mathbf{z}\dvtx 1-\varepsilon \leq q(\mathbf
{z})\leq1+\varepsilon \}\bigr)\\
&=&\nu(S)\bigl((1+\varepsilon )^{\alpha /(d+\alpha )}-(1-\varepsilon
)^{\alpha /(d+\alpha )}\bigr).
\end{eqnarray*}
Hence,
$\lim_{n\rightarrow\infty }\frac{P(\bar{Q}_n\bigtriangleup
\widetilde{Q}_n)}{p}=0$.
\end{pf}

The following proposition shows uniform consistency of $\widehat{\psi
}_n$ and
might be of independent interest. There is an abundant literature on
density estimation for directional data. In particular, uniform
consistency of density estimators for directional data has been
established in \citet{BaiRaoZha88}. Here, however, the data do
not have a \textit{fixed} probability density on $\Theta$: the
density $\psi$ is defined via a limit relation. Hence, $\psi$ is only
an approximate model for the directional data. As a consequence, a more
general result is required.

\begin{prop}\label{prop2}
As $n\rightarrow\infty $,
\[
\sup_{\mathbf{w}\in\Theta} |\widehat{\psi}_n(\mathbf{w})-\psi
(\mathbf{w})
|\stackrel{\mathbb{P}}{\rightarrow}0.
\]
\end{prop}

\begin{pf}
It is easy to see that, for any $\eta>0$, there exists a function
\[
K^*=\sum_{j=1}^m \alpha _j 1_{[r_{j-1}, r_j)}
\]
with $1 \geq\alpha _1\geq\alpha _2\geq\cdots\geq\alpha _m\geq
0$ and
$0=r_0<r_1<\cdots<r_m=1$, such that
\[
\sup_{u\in[0,1]}|K(u)-K^*(u)|\leq\eta.
\]

Write $U_i=1-H(R_i)$, $i=1,\ldots,n$, and denote the corresponding order
statistics with $U_{i: n}$. Let $\widetilde{P}$ be the probability
measure on
$\Theta
\times(0,1)$ corresponding to $(\mathbf{W}_1, U_1)$ and let
$\widetilde{P}_n$ be the
empirical
measure of the $(\mathbf{W}_i, U_i)$ $i=1,\ldots, n$.
Define
\[
\psi_n^*(\mathbf{w})=\frac{c(h, K)}{k}\sum_{i=1}^n K^*\biggl(\frac
{1-\mathbf{w}^T\mathbf{W}_i}{h}\biggr)1_{[R_i>R_{n-k: n}]}
\]
and
\[
\psi_{n,j}^*(\mathbf{w})=\frac{nc(h, K)}{k}\widetilde
{P}_n\bigl(D_{\mathbf{w},j}\times
(0,U_{k: n}]\bigr)\vadjust{\goodbreak}
\]
with $D_{\mathbf{w},j}= \{\mathbf{v}\in\Theta\dvtx 1-hr_{j}<\mathbf{w}^T\mathbf{v}
\leq1-hr_{j-1} \}$.
Observe that $\psi_n^*(\mathbf{w})=\sum_{j=1}^m\alpha _j\psi
_{n,j}^*(\mathbf{w})$. Also write
\[
\psi_{n,j}(\mathbf{w})=\frac{nc(h, K)}{k}\widetilde{P}\bigl(D_{\mathbf{w},j}\times(0,U_{k: n}]\bigr).
\]

Let $\varepsilon >0$. It is sufficient to show that for large $n$
\begin{eqnarray}
\mathbb{P}\Biggl(\sup_{\mathbf{w}\in\Theta} \Biggl|\widehat{\psi
}_n(\mathbf{w})-\sum
_{j=1}^m\alpha _j\psi_{n,j}(\mathbf{w}) \Biggr|\geq2\varepsilon \Biggr)&\leq&
2\varepsilon \label{(p2.1)},\\
\mathbb{P}\Biggl(\sup_{\mathbf{w}\in\Theta} \Biggl|\sum_{j=1}^m\alpha_j\bigl(\psi
_{n,j}(\mathbf{w})-c(h, K)\Psi(D_{\mathbf{w},j})\bigr) \Biggr|\geq
2\varepsilon \Biggr)&\leq&\varepsilon, \label{(p2.2)}\\
\sup_{\mathbf{w}\in\Theta} \Biggl|c(h, K)\sum_{j=1}^m\alpha _j\Psi
(D_{\mathbf{w},j})-\psi(\mathbf{w}) \Biggr|&\leq&\varepsilon . \label{(p2.3)}
\end{eqnarray}
For $\mathbf{w}\in\Theta$ and $\delta \in(0,1)$, write
$\mathcal{C}_{\delta }=\{C_\mathbf{w}(a)\dvtx \mathbf{w}\in\Theta
,a\leq\delta
\}$. Note that, as $n\rightarrow\infty $,
\begin{equation}\label{(p2.4)}
\sup_{C\in\mathcal{C}_{1}, 0<s\leq2}\frac{1}{\lambda (C)}
\biggl|\frac{n}{k}\widetilde{P}\bigl(C\times(0, sk/n ]\bigr)-s\Psi(C) \biggr|\rightarrow
0.
\end{equation}
This readily follows from
\begin{eqnarray*}
\frac{n}{k}\widetilde{P}\bigl(C\times(0, sk/n ]\bigr)
&=&\frac{n}{k}\mathbb{P}\biggl(\mathbf{W}\in C, R\geq U\biggl(\frac{n}{sk}\biggr)\biggr)\\
&=&\frac{n}{k}\int_{U({n}/{(sk)})}^{\infty }\int_C\frac
{f(r\mathbf{w})}{r^{-d}V(r)}\,d\lambda (\mathbf{w})\,r^{-1}V(r)\,dr
\end{eqnarray*}
and (\ref{(4.1)}) and (\ref{(4.5)}).

Now we prove (\ref{(p2.1)}).
It is easy to show that
\[
c(h,K)=\biggl(\frac{2\pi^{{(d-1)}/{2}}}{\Gamma({(d-1)}/{2})}
\int_{1-h}^1K\biggl(\frac{1-t}{h}\biggr)(1-t^2)^{{(d-3)}/{2}}\,dt\biggr)^{-1}
\]
%
and hence
\begin{eqnarray}
\limsup_{h\downarrow0}c(h, K)\lambda (C_\mathbf{w}(h))<\infty .
\label{(p2.0)}
\end{eqnarray}
We have
\begin{eqnarray} \label{(p2.1.2)}
&&|\widehat{\psi}_n(\mathbf{w})-\psi_n^*(\mathbf{w}) |\nonumber \\
&&\qquad =\frac{c(h, K)}{k} \Biggl|\sum_{i=1}^n\biggl(K\biggl(\frac{1-\mathbf{w}^T\mathbf{W}
_i}{h}\biggr)-K^*\biggl(\frac{1-\mathbf{w}^T\mathbf
{W}_i}{h}\biggr)\biggr)1_{[R_i>R_{n-k: n}]} \Biggr|
\nonumber\\
&&\qquad \leq \frac{c(h, K)}{k}\sum_{i=1}^n\eta1_{[\mathbf{W}_i\in
C_\mathbf{w}(h), R_i>R_{n-k: n}]} \\
&&\qquad \leq \eta\frac{n c(h, K)}{k}\widetilde{P}\bigl(C_\mathbf{w}(h)\times
(0,U_{k: n}]\bigr) \nonumber\\
&&\qquad \quad {} +\eta\frac{n c(h, K)}{k} \bigl|(\widetilde{P}_n-\widetilde
{P})\bigl(C_\mathbf{w}(h)\times
(0,U_{k: n}]\bigr) \bigr|.\nonumber
\end{eqnarray}
By (\ref{(p2.4)}), for $\eta$ small enough the first term is
less than $\varepsilon $, with probability tending to one, uniformly in
$\mathbf{w}\in\Theta$. Also,
\begin{eqnarray} \label{(p2.1.1)}
&&\Biggl|\psi_n^*(\mathbf{w})-\sum_{j=1}^m\alpha _j\psi_{n,j}(\mathbf{w})
\Biggr|\nonumber\\
&&\qquad \leq \sum_{j=1}^m\alpha _j |\psi_{n,j}^*(\mathbf{w})-\psi
_{n,j}(\mathbf{w}) | \\
&&\qquad \leq \sum_{j=1}^m\alpha _j \frac{nc(h,K)}{k} \bigl|(\widetilde
{P}_n-\widetilde{P})\bigl(D_{\mathbf{w},j}\times(0,U_{k: n}]\bigr) \bigr|.\nonumber
\end{eqnarray}
From (\ref{(p2.1.1)}), (\ref{(p2.1.2)}) and (\ref{(p2.0)}), we see
that for a proof of
(\ref{(p2.1)}) it remains to show that
\[
\frac{n}{k\lambda (C_\mathbf{w}(h))}\sup_{\mathbf{w}\in\Theta
}\sup
_{0<a\leq1}
\bigl|(\widetilde{P}_n-\widetilde{P})\bigl(C_\mathbf{w}(ah)\times(0,U_{k: n}]\bigr)
\bigr|\stackrel{\mathbb{P}}{\rightarrow}0.
\]
It can be shown that there exists a constant $c=c(d)$ and
finitely many $\mathbf{w}_l$, $l=1,\ldots,l_h$ such that $l_h=O(c(h, K))$
as $h\downarrow0$, and for every $\mathbf{w}\in\Theta$ and $0<a\leq1$
\[
C_\mathbf{w}(ah)\in C_{\mathbf{w}_l}(ch)\qquad  \mbox{for some } l.
\]
Hence for $\varepsilon _1>0$,
\begin{eqnarray*}
&&\mathbb{P}\biggl(\frac{n}{k\lambda (C_\mathbf{w}(h))}\sup_{\mathbf{w}\in\Theta
}\sup_{0<a\leq1}
\bigl|(\widetilde{P}_n-\widetilde{P})\bigl(C_\mathbf{w}(ah)\times(0,U_{k: n}]\bigr)
\bigr|\geq\varepsilon _1\biggr)\\
&&\qquad \leq\mathbb{P}\Bigl(\max_{1\leq l\leq l_h}\mathop{\sup_{C\subset
C_{\mathbf{w}_l}(ch)}}_{ C\in\mathcal{C}_h}\sup_{0<s\leq2}
\bigl|(\widetilde{P}_n-\widetilde{P})\bigl(C\times(0,sk/n]\bigr) \bigr|\geq\varepsilon
_1k/n\lambda (C_\mathbf{w}(h))\Bigr)\\
&&\qquad \quad {}+\mathbb{P}(U_{k: n}>2k/n)\\
&&\qquad \leq\sum_{l=1}^{l_h}\mathbb{P}\Bigl(\mathop{\sup_{C\subset
C_{\mathbf{w}_l}(ch)}}_{ C\in\mathcal{C}_h}\sup_{0<s\leq2}
\bigl|(\widetilde{P}_n-\widetilde{P})\bigl(C\times(0,sk/n]\bigr) \bigr|\geq\varepsilon
_1k/n\lambda (C_\mathbf{w}(h))\Bigr)\\
&&\qquad \quad {}+\mathbb{P}(U_{k: n}>2k/n).
\end{eqnarray*}
The latter probability tends to $0$, so it suffices to consider
the sum of the $l_h$ probabilities. Write $b=\varepsilon _1k\lambda
(C_\mathbf{w}(h))$. Fix $l$ and define
$N=n\widetilde{P}_n(C_{\mathbf{w}_l}(ch)\times(0,2k/n])$,
$\mu=n\widetilde{P}(C_{\mathbf{w}_l}(ch)\times(0,2k/n])$. Define the
conditional probability measure
$\widetilde{P}_c=\frac{n \widetilde{P}}{\mu}$ on $C_{\mathbf{w}_l}(ch)\times(0,2k/n]$
and let $\widetilde{P}_{c,r}$ be the corresponding empirical measure,
based on
$r$ observations. We have
\begin{eqnarray} \label{(p2.1.3)}
&&\mathbb{P}\Bigl(\mathop{\sup_{C\subset C_{\mathbf{w}_l}(ch)}}_{ C\in
\mathcal{C}_h}\sup_{0<s\leq2}
n \bigl|(\widetilde{P}_n-\widetilde{P})\bigl(C\times(0,sk/n]\bigr) \bigr|\geq b\Bigr)
\nonumber\\[-1pt]
&&\qquad \leq\sum_{r=\lceil\mu-b/3\rceil}^{r=\lfloor\mu+b/3\rfloor}
\mathbb{P}\Bigl(\mathop{\sup_{C\subset C_{\mathbf{w}_l}(ch)}}_{ C\in
\mathcal{C}_h}\sup_{0<s\leq2}
n \bigl|(\widetilde{P}_n-\widetilde{P})\bigl(C\times(0,sk/n]\bigr) \bigr|\geq b
\big|N=r\Bigr)\nonumber\\[-1pt]
&&\qquad \quad\hphantom{\sum_{r=\lceil\mu-b/3\rceil}^{r=\lfloor\mu+b/3\rfloor}} {}\times\mathbb{P}(N=r) +\mathbb{P}(|N-\mu|\geq b/3) \nonumber\\[-1pt]
&&\qquad \leq\sum_{r=\lceil\mu-b/3\rceil}^{r=\lfloor\mu+b/3\rfloor}
\mathbb{P}\biggl(\mathop{\sup_{C\subset C_{\mathbf{w}_l}(ch)}}_{ C\in
\mathcal{C}_h}\sup_{0<s\leq2}
n \biggl|\biggl(\widetilde{P}_n-\frac{N}{\mu}\widetilde{P}\biggr)\bigl(C\times(0,sk/n]\bigr)
\biggr|\nonumber \\[-1pt]
&&\hspace*{252pt}\geq\frac{b}{2}
\Big|N=r\biggr)\mathbb{P}(N=r)\nonumber \\[-1pt]
&&\qquad \quad {} +\sum_{r=\lceil\mu-b/3\rceil}^{r=\lfloor\mu+b/3\rfloor}
\mathbb{P}\biggl(\mathop{\sup_{C\subset C_{\mathbf{w}_l}(ch)}}_{ C\in
\mathcal{C}_h}\sup_{0<s\leq2}
n \biggl|\frac{(N-\mu)}{\mu}\widetilde{P}\bigl(C\times(0,sk/n]\bigr) \biggr|\nonumber \\[-1pt]
&&\hspace*{255pt}\geq\frac{b}{2}
\Big|N=r\biggr)\mathbb{P}(N=r)\nonumber \\[-1pt]
&&\qquad \quad {}+\mathbb{P}(|N-\mu|\geq b/3) \\[-1pt]
&&\qquad \leq \sum_{r=\lceil\mu-b/3\rceil}^{r=\lfloor\mu+b/3\rfloor}
\mathbb{P}\biggl(\mathop{\sup_{C\subset C_{\mathbf{w}_l}(ch)}}_{ C\in
\mathcal{C}_h} \sup_{0<s\leq2}
r \bigl|(\widetilde{P}_{c,r}-\widetilde{P}_c)\bigl(C\times(0,sk/n]\bigr) \bigr|\geq
\frac{b}{2}\biggr)\nonumber\\[-1pt]
&&\qquad \quad\hphantom{\sum_{r=\lceil\mu-b/3\rceil}^{r=\lfloor\mu+b/3\rfloor}} {}\times\mathbb{P}
(N=r)\nonumber\\[-1pt]
&&\qquad \quad {} + \sum_{r=\lceil\mu-b/3\rceil}^{r=\lfloor\mu+b/3\rfloor}
\mathbb{P}\biggl(|r-\mu|\geq\frac{b}{2}\biggr)\mathbb{P}(N=r)+\mathbb
{P}(|N-\mu|\geq b/3).\nonumber
\end{eqnarray}
Note that the first probability of the second sum in the
right
side of (\ref{(p2.1.3)}) is equal to~0.
From Bennett's inequality [cf.\ \citet{ShoWel86}, page~851],
it follows that for some constant $c_1$
\[
\mathbb{P}(|N-\mu|\geq b/3)\leq2\exp\biggl(-\varepsilon _1^2c_1\frac
{k}{c(h, K)}\biggr).\vadjust{\goodbreak}
\]
Hence, since $l_h=O(c(h, K))$,
\[
\sum_{l=1}^{l_h}\mathbb{P}(|N-\mu|\geq b/3)
=O\biggl(c(h, K)\exp\biggl(-\varepsilon _1^2c_1\frac{k}{c(h, K)}\biggr)\biggr)=o(1).
\]
To complete the proof of (\ref{(p2.1)}), we need to consider
the first sum in the right side of (\ref{(p2.1.3)}). For the
first probability in there, we use Corollary 2.9 in \citet
{Ale84}, a
good probability bound for empirical processes on VC classes.
We obtain as an upper bound
\[
16\exp\biggl(-\frac{b^2}{4r}\biggr).
\]
Using $r\leq\mu+b/3$, we find for some constant $c_2$
\begin{eqnarray*}
&&\sum_{l=1}^{l_h}\sum_{r=\lceil\mu-b/3\rceil}^{r=\lfloor\mu+b/3\rfloor}
\mathbb{P}\biggl(\mathop{\sup_{C\subset C_{\mathbf{w}_l}(ch)}}_{ C\in
\mathcal
{C}_h}\sup_{0<s\leq2}
r \bigl|(\widetilde{P}_{c,r}-\widetilde{P}_c)\bigl(C\times(0,sk/n]\bigr) \bigr|\geq
\frac{b}{2}\biggr)\\
&&\hphantom{\sum_{l=1}^{l_h}\sum_{r=\lceil\mu-b/3\rceil}^{r=\lfloor\mu+b/3\rfloor}}
{}\times\mathbb{P}(N=r)\\
&&\qquad \leq16 \sum_{l=1}^{l_h}\sum_{r=\lceil\mu-b/3\rceil}^{r=\lfloor
\mu+b/3\rfloor}
\exp\biggl(-\varepsilon _1^2c_2\frac{k}{c(h, K)}\biggr)\mathbb{P}(N=r)\\
&&\qquad \leq16 \sum_{l=1}^{l_h} \exp\biggl(-\varepsilon _1^2c_2\frac{k}{c(h, K)}\biggr)\\
&&\qquad =o(1).
\end{eqnarray*}

Next, we show (\ref{(p2.2)}). From (\ref{(p2.0)}) and
(\ref{(p2.4)}), we obtain for $\varepsilon _2>0$ small enough,
\begin{eqnarray*}
&&\sup_{\mathbf{w}\in\Theta} \Biggl|\sum_{j=1}^m\alpha _j\bigl(\psi
_{n,j}(\mathbf{w})-c(h, K)\Psi(D_{\mathbf{w},j})\bigr) \Biggr|\\
&&\qquad =\sup_{\mathbf{w}\in\Theta} \Biggl|\sum_{j=1}^m\alpha _j c(h,
K)\bigl(n/k\widetilde{P}
\bigl(D_{\mathbf{w},j}\times(0,U_{k: n}]\bigr)-\Psi(D_{\mathbf{w},j})\bigr) \Biggr|\\
&&\qquad \leq \varepsilon _2\sum_{j=1}^m\alpha _j c(h, K)\lambda
(C_\mathbf{w}(h))+
\sup_{\mathbf{w}\in\Theta} \Biggl|\sum_{j=1}^m\alpha _jc(h,
K)(nU_{k: n}/k-1)\Psi(D_{\mathbf{w},j}) \Biggr|\\
&&\qquad \leq\varepsilon + \biggl|\frac{n}{k}U_{k: n}-1 \biggr|\sum_{j=1}^m\alpha
_jc(h, K)\lambda
(C_\mathbf{w}(h))\sup_{\mathbf{w}\in\Theta}\psi(\mathbf{w})<2\varepsilon
\end{eqnarray*}
with probability tending to one.

It remains to prove (\ref{(p2.3)}). It is readily seen that
$\int_{C_\mathbf{w}(h)}K^*(\frac{1-\mathbf{w}^T\mathbf
{v}}{h})\,d\lambda (\mathbf{v}
)=\sum_{j=1}^m\alpha _j\lambda (D_{\mathbf{w},j})$.
Hence, for $\varepsilon _3>0$ small enough
\begin{eqnarray*}
&&\sup_{\mathbf{w}\in\Theta} \Biggl|c(h, K)\sum_{j=1}^m\alpha _j\Psi
(D_{\mathbf{w},j})-\psi(\mathbf{w}) \Biggr|\\[-1pt]
&&\qquad \leq \sup_{\mathbf{w}\in\Theta}\psi(\mathbf{w}) \Biggl|c(h, K)\sum
_{j=1}^m\alpha _j\lambda (D_{\mathbf{w},j})-1 \Biggr|+\varepsilon _3 c(h,
K)\sum_{j=1}^m\alpha
_j\lambda (D_{\mathbf{w},j})\\[-1pt]
&&\qquad \leq \sup_{\mathbf{w}\in\Theta}\psi(\mathbf{w}) \biggl|\frac{\int
_{C_\mathbf{w}(h)}K^*({(1-\mathbf{w}^T\mathbf{v})}/{h})\,d\lambda
(\mathbf{v})}
{\int_{C_\mathbf{w}(h)}K({(1-\mathbf{w}^T\mathbf
{v})}/{h})\,d\lambda (\mathbf{v})}-1
\biggr|\\[-1pt]
&&\qquad \quad {}+\varepsilon _3 c(h, K)\lambda (C_{\mathbf{w}}(h))\sum
_{j=1}^m\alpha _j\\[-1pt]
&&\qquad \leq \eta c(h, K)\lambda (C_{\mathbf{w}}(h))\sup_{\mathbf{w}\in
\Theta
}\psi(\mathbf{w})+\varepsilon _3 c(h, K)\lambda (C_{\mathbf{w}}(h))\sum_{j=1}^m\alpha
_j\\[-1pt]
&&\qquad \leq \varepsilon .
\end{eqnarray*}
\upqed
\end{pf}

From Proposition~\ref{prop2} and the consistency of $\widehat\alpha$, we obtain
immediately, as $n\rightarrow\infty $,
\[
\widehat{\nu(S)}\stackrel{\mathbb{P}}{\rightarrow}\nu(S)
\]
and, for $\varepsilon>0$,
\begin{equation}\label{inclu}
\mathbb{P}\bigl( (1+\varepsilon )S\subset\widehat{S}\subset
(1-\varepsilon )S\bigr)\to1.
\end{equation}

\begin{prop}\label{prop3}
As $n\rightarrow\infty $,
\[
\frac{P(\widetilde{Q}_n\triangle\widehat{Q}_n)}{p} \stackrel{\mathbb{P}}{\rightarrow} 0.
\]
\end{prop}

\begin{pf}
Note that as $n\rightarrow\infty $, we have
\begin{eqnarray*}
\widehat{U}\biggl(\frac{n}{k}\biggr)\Big/U\biggl(\frac{n}{k}\biggr)&\stackrel{\mathbb{P}}{\rightarrow}&1,\\[-1pt]
(\widehat{\nu(S)})^{{1}/{\widehat{\alpha
}}}&\stackrel{\mathbb{P}}{\rightarrow}&
(\nu(S))^{1/\alpha },
\\[-1pt]
\biggl(\frac{k}{np}\biggr)^{1/\widehat{\alpha }-1/\alpha }
&=&\exp\biggl(\frac{\sqrt{k}(\alpha -\widehat{\alpha })}{\widehat
{\alpha }\alpha }
\biggl(\frac{\log k}{\sqrt{k}}-\frac{\log(np)}{\sqrt{k}}\biggr)\biggr)
\stackrel{\mathbb{P}}{\rightarrow}1.
\end{eqnarray*}
Combining these three limit relations, we obtain
\[
\frac{\widehat{U}({n}/{k})({k\widehat{\nu
(S)}}/{(np)})^{{1}/{\widehat{\alpha
}}}}{U({n}/{k})({k\nu(S)}/{(np)})^{{1}/{\alpha }}}
\stackrel{\mathbb{P}}{\rightarrow} 1.\vadjust{\goodbreak}
\]
This and (\ref{inclu}) yields
that with probability tending to one, as $n \to\infty $,
\[
(1+\varepsilon )^2\widetilde{Q}_n\subset\widehat{Q}_n\subset
(1-\varepsilon )^2\widetilde{Q}_n.
\]
Then,
\[
\frac{P(\widetilde{Q}_n\triangle\widehat{Q}_n)}{p}
\leq\frac{1}{p}P\biggl(U\biggl(\frac{n}{k}\biggr)\biggl(\frac{k\nu(S)}{np}\biggr)^{
{1}/{\alpha
}}\bigl((1-\varepsilon )^2S\setminus(1+\varepsilon )^2S\bigr)\biggr),
\]
and, by (\ref{(4.32)}), the latter expression tends to
\begin{eqnarray*}
&&\nu\bigl((1-\varepsilon )^2S \setminus(1+\varepsilon
)^2S\bigr)/\nu(S)\\
&&\qquad =\nu\bigl((1-\varepsilon )^2S\bigr)/\nu(S) -\nu\bigl((1+\varepsilon )^2S\bigr)/\nu
(S)\\
&&\qquad =(1-\varepsilon )^{-2\alpha }-(1+\varepsilon )^{-2\alpha },
\end{eqnarray*}
which in
turn tends to 0, as $\varepsilon \downarrow0$.
\end{pf}

\begin{pf*}{Proof of Theorem~\ref{thm1}}
The result follows from
Propositions~\ref{prop1} and~\ref{prop3}.~%
\end{pf*}

\section*{Acknowledgments}
We thank two referees for many insightful comments that led to an
improved version of the paper. We are grateful to Kees Koedijk, Roger
Laeven, Ronald Mahieu and Chen Zhou for discussions of the financial
application.


%

\printaddresses


\begin{thebibliography}{20}

\bibitem[\protect\citeauthoryear{Alexander}{1984}]{Ale84}
%
\begin{barticle}[mr]
\bauthor{\bsnm{Alexander},~\bfnm{Kenneth~S.}\binits{K.~S.}}
(\byear{1984}).
\btitle{Probability inequalities for empirical processes and a law of the
iterated logarithm}.
\bjournal{Ann. Probab.}
\bvolume{12}
\bpages{1041--1067}.
\bnote{[Corr.: \textbf{15} (1987) 428--430.]}
\bid{issn={0091-1798}, mr={0757769}}
\bptnote{check related}%
\end{barticle}
%
\endbibitem

\bibitem[\protect\citeauthoryear{Bai, Rao and Zhao}{1988}]{BaiRaoZha88}
%
\begin{barticle}[mr]
\bauthor{\bsnm{Bai},~\bfnm{Z.~D.}\binits{Z.~D.}},
\bauthor{\bsnm{Rao},~\bfnm{C.~Radhakrishna}\binits{C.~R.}} \AND
\bauthor{\bsnm{Zhao},~\bfnm{L.~C.}\binits{L.~C.}}
(\byear{1988}).
\btitle{Kernel estimators of density function of directional data}.
\bjournal{J. Multivariate Anal.}
\bvolume{27}
\bpages{24--39}.
\bid{doi={10.1016/0047-259X(88)90113-3}, issn={0047-259X}, mr={0971170}}
\end{barticle}
%
\endbibitem

\bibitem[\protect\citeauthoryear{Ba{\'{\i}}llo, Cuesta-Albertos and
Cuevas}{2001}]{BalCueCue01}
%
\begin{barticle}[mr]
\bauthor{\bsnm{Ba{\'{\i}}llo},~\bfnm{Amparo}\binits{A.}},
\bauthor{\bsnm{Cuesta-Albertos},~\bfnm{Juan~A.}\binits{J.~A.}} \AND
\bauthor{\bsnm{Cuevas},~\bfnm{Antonio}\binits{A.}}
(\byear{2001}).
\btitle{Convergence rates in nonparametric estimation of level sets}.
\bjournal{Statist. Probab. Lett.}
\bvolume{53}
\bpages{27--35}.
\bid{doi={10.1016/S0167-7152(01)00006-2}, issn={0167-7152}, mr={1843338}}
\end{barticle}
%
\endbibitem

\bibitem[\protect\citeauthoryear{Danielsson et~al.}{2001}]{Danetal01}
%
\begin{barticle}[mr]
\bauthor{\bsnm{Danielsson},~\bfnm{J.}\binits{J.}}, \bauthor
{\bparticle{de}
\bsnm{Haan},~\bfnm{L.}\binits{L.}},
\bauthor{\bsnm{Peng},~\bfnm{L.}\binits{L.}} \AND\bauthor
{\bparticle{de}
\bsnm{Vries},~\bfnm{C.~G.}\binits{C.~G.}}
(\byear{2001}).
\btitle{Using a bootstrap method to choose the sample fraction in tail index
estimation}.
\bjournal{J. Multivariate Anal.}
\bvolume{76}
\bpages{226--248}.
\bid{doi={10.1006/jmva.2000.1903}, issn={0047-259X}, mr={1821820}}
\end{barticle}
%
\endbibitem

\bibitem[\protect\citeauthoryear{de~Haan and Ferreira}{2006}]{deHFer06}
%
\begin{bbook}[mr]
\bauthor{\bparticle{de} \bsnm{Haan},~\bfnm{Laurens}\binits{L.}}
\AND
\bauthor{\bsnm{Ferreira},~\bfnm{Ana}\binits{A.}}
(\byear{2006}).
\btitle{Extreme Value Theory: An Introduction}.
\bpublisher{Springer}, \baddress{New York}.
\bid{mr={2234156}}
\end{bbook}
%
\endbibitem

\bibitem[\protect\citeauthoryear{de~Haan and Resnick}{1987}]{deHRes87}
%
\begin{barticle}[mr]
\bauthor{\bparticle{de} \bsnm{Haan},~\bfnm{L.}\binits{L.}} \AND
\bauthor{\bsnm{Resnick},~\bfnm{S.}\binits{S.}}
(\byear{1987}).
\btitle{On regular variation of probability densities}.
\bjournal{Stochastic Process. Appl.}
\bvolume{25}
\bpages{83--93}.
\bid{doi={10.1016/0304-4149(87)90191-8}, issn={0304-4149}, mr={0904266}}
\end{barticle}
%
\endbibitem

\bibitem[\protect\citeauthoryear{Dekkers, Einmahl and
de~Haan}{1989}]{DekEindeH89}
%
\begin{barticle}[mr]
\bauthor{\bsnm{Dekkers},~\bfnm{A.~L.~M.}\binits{A.~L.~M.}},
\bauthor{\bsnm{Einmahl},~\bfnm{J.~H.~J.}\binits{J.~H.~J.}} \AND
\bauthor{\bparticle{de} \bsnm{Haan},~\bfnm{L.}\binits{L.}}
(\byear{1989}).
\btitle{A moment estimator for the index of an extreme-value distribution}.
\bjournal{Ann. Statist.}
\bvolume{17}
\bpages{1833--1855}.
\bid{doi={10.1214/aos/1176347397}, issn={0090-5364}, mr={1026315}}
\end{barticle}
%
\endbibitem

\bibitem[\protect\citeauthoryear{Einmahl, Li and Liu}{2009}]{EinLiLiu09}
%
\begin{barticle}[mr]
\bauthor{\bsnm{Einmahl},~\bfnm{John H.~J.}\binits{J.~H.~J.}},
\bauthor{\bsnm{Li},~\bfnm{Jun}\binits{J.}} \AND
\bauthor{\bsnm{Liu},~\bfnm{Regina~Y.}\binits{R.~Y.}}
(\byear{2009}).
\btitle{Thresholding events of extreme in simultaneous monitoring of multiple
risks}.
\bjournal{J. Amer. Statist. Assoc.}
\bvolume{104}
\bpages{982--992}.
\bid{doi={10.1198/jasa.2009.ap08329}, issn={0162-1459}, mr={2562001}}
\end{barticle}
%
\endbibitem

\bibitem[\protect\citeauthoryear{Hall, Watson and
Cabrera}{1987}]{HalWatCab87}
%
\begin{barticle}[mr]
\bauthor{\bsnm{Hall},~\bfnm{Peter}\binits{P.}},
\bauthor{\bsnm{Watson},~\bfnm{G.~S.}\binits{G.~S.}} \AND
\bauthor{\bsnm{Cabrera},~\bfnm{Javier}\binits{J.}}
(\byear{1987}).
\btitle{Kernel density estimation with spherical data}.
\bjournal{Biometrika}
\bvolume{74}
\bpages{751--762}.
\bid{doi={10.1093/biomet/74.4.751}, issn={0006-3444}, mr={0919843}}
\end{barticle}
%
\endbibitem

\bibitem[\protect\citeauthoryear{Hashorva}{2006}]{Has06}
%
\begin{barticle}[mr]
\bauthor{\bsnm{Hashorva},~\bfnm{Enkelejd}\binits{E.}}
(\byear{2006}).
\btitle{On the regular variation of elliptical random vectors}.
\bjournal{Statist. Probab. Lett.}
\bvolume{76}
\bpages{1427--1434}.
\bid{doi={10.1016/j.spl.2006.02.014}, issn={0167-7152}, mr={2245561}}
\end{barticle}
%
\endbibitem

\bibitem[\protect\citeauthoryear{Hill}{1975}]{Hil75}
%
\begin{barticle}[mr]
\bauthor{\bsnm{Hill},~\bfnm{Bruce~M.}\binits{B.~M.}}
(\byear{1975}).
\btitle{A simple general approach to inference about the tail of a
distribution}.
\bjournal{Ann. Statist.}
\bvolume{3}
\bpages{1163--1174}.
\bid{issn={0090-5364}, mr={0378204}}
\end{barticle}\vadjust{\goodbreak}
%
\endbibitem

\bibitem[\protect\citeauthoryear{Jessen and Mikosch}{2006}]{JesMik06}
%
\begin{barticle}[mr]
\bauthor{\bsnm{Jessen},~\bfnm{Anders~Hedegaard}\binits{A.~H.}} \AND
\bauthor{\bsnm{Mikosch},~\bfnm{Thomas}\binits{T.}}
(\byear{2006}).
\btitle{Regularly varying functions}.
\bjournal{Publ. Inst. Math. (Beograd) (N.S.)}
\bvolume{80(94)}
\bpages{171--192}.
\bid{doi={10.2298/PIM0694171J}, issn={0350-1302}, mr={2281913}}
\end{barticle}
%
\endbibitem

\bibitem[\protect\citeauthoryear{M{\"u}ller and Sawitzki}{1991}]{MulSaw91}
%
\begin{barticle}[mr]
\bauthor{\bsnm{M{\"u}ller},~\bfnm{D.~W.}\binits{D.~W.}} \AND
\bauthor{\bsnm{Sawitzki},~\bfnm{G.}\binits{G.}}
(\byear{1991}).
\btitle{Excess mass estimates and tests for multimodality}.
\bjournal{J.~Amer. Statist. Assoc.}
\bvolume{86}
\bpages{738--746}.
\bid{issn={0162-1459}, mr={1147099}}
\end{barticle}
%
\endbibitem

\bibitem[\protect\citeauthoryear{Polonik}{1995}]{Pol95}
%
\begin{barticle}[mr]
\bauthor{\bsnm{Polonik},~\bfnm{Wolfgang}\binits{W.}}
(\byear{1995}).
\btitle{Measuring mass concentrations and estimating density contour
clusters---an excess mass approach}.
\bjournal{Ann. Statist.}
\bvolume{23}
\bpages{855--881}.
\bid{doi={10.1214/aos/1176324626}, issn={0090-5364}, mr={1345204}}
\end{barticle}
%
\endbibitem

\bibitem[\protect\citeauthoryear{Rigollet and Vert}{2009}]{RigVer09}
%
\begin{barticle}[mr]
\bauthor{\bsnm{Rigollet},~\bfnm{Philippe}\binits{P.}} \AND
\bauthor{\bsnm{Vert},~\bfnm{R{\'e}gis}\binits{R.}}
(\byear{2009}).
\btitle{Optimal rates for plug-in estimators of density level sets}.
\bjournal{Bernoulli}
\bvolume{15}
\bpages{1154--1178}.
\bid{doi={10.3150/09-BEJ184}, issn={1350-7265}, mr={2597587}}
\end{barticle}
%
\endbibitem

\bibitem[\protect\citeauthoryear{Rva{\v{c}}eva}{1962}]{Rva62}
%
\begin{barticle}[mr]
\bauthor{\bsnm{Rva{\v{c}}eva},~\bfnm{E.~L.}\binits{E.~L.}}
(\byear{1962}).
\btitle{On domains of attraction of multi-dimensional distributions}.
\bjournal{Select. {T}ransl. {M}ath. {S}tatist. {P}robab.}
\bvolume{2}
\bpages{183--205}.
\bid{mr={0150795}}
\end{barticle}
%
\endbibitem

\bibitem[\protect\citeauthoryear{Shorack and Wellner}{1986}]{ShoWel86}
%
\begin{bbook}[mr]
\bauthor{\bsnm{Shorack},~\bfnm{Galen~R.}\binits{G.~R.}} \AND
\bauthor{\bsnm{Wellner},~\bfnm{Jon~A.}\binits{J.~A.}}
(\byear{1986}).
\btitle{Empirical Processes with Applications to Statistics}.
\bpublisher{Wiley}, \baddress{New York}.
\bid{mr={0838963}}
\end{bbook}
%
\endbibitem

\bibitem[\protect\citeauthoryear{Smith}{1987}]{Smi87}
%
\begin{barticle}[mr]
\bauthor{\bsnm{Smith},~\bfnm{Richard~L.}\binits{R.~L.}}
(\byear{1987}).
\btitle{Estimating tails of probability distributions}.
\bjournal{Ann. Statist.}
\bvolume{15}
\bpages{1174--1207}.
\bid{doi={10.1214/aos/1176350499}, issn={0090-5364}, mr={0902252}}
\end{barticle}
%
\endbibitem

\bibitem[\protect\citeauthoryear{Tsybakov}{1997}]{Tsy97}
%
\begin{barticle}[mr]
\bauthor{\bsnm{Tsybakov},~\bfnm{A.~B.}\binits{A.~B.}}
(\byear{1997}).
\btitle{On nonparametric estimation of density level sets}.
\bjournal{Ann. Statist.}
\bvolume{25}
\bpages{948--969}.
\bid{doi={10.1214/aos/1069362732}, issn={0090-5364}, mr={1447735}}
\end{barticle}
%
\endbibitem

\end{thebibliography}
\end{document}